\documentclass[12pt]{article}
\usepackage{}

\usepackage{epsfig}
\usepackage{latexsym}
\usepackage{caption}
\usepackage{amsfonts}
\usepackage{amssymb}
\usepackage{mathrsfs}
\usepackage{amsmath}
\usepackage{enumerate}
\usepackage{graphics}
\usepackage{MnSymbol}
\usepackage{float}
\usepackage{pict2e}
\usepackage{subfigure}

\usepackage{color}

\makeatletter

\newcommand{\Rmnum}[1]{\expandafter\@slowromancap\romannumeral #1@}

\makeatother

\DeclareMathOperator{\ch}{ch}
\DeclareMathOperator{\AT}{AT}
\DeclareMathOperator{\p}{P}

\begin{document}

\newtheorem{theorem}{Theorem}
\newtheorem{observation}[theorem]{Observation}
\newtheorem{corollary}[theorem]{Corollary}
\newtheorem{algorithm}[theorem]{Algorithm}
\newtheorem{definition}[theorem]{Definition}
\newtheorem{guess}[theorem]{Conjecture}
\newtheorem{claim}[theorem]{Claim}
\newtheorem{problem}[theorem]{Problem}
\newtheorem{question}[theorem]{Question}
\newtheorem{lemma}[theorem]{Lemma}
\newtheorem{proposition}[theorem]{Proposition}
\newtheorem{fact}[theorem]{Fact}

\makeatletter
  \newcommand\figcaption{\def\@captype{figure}\caption}
  \newcommand\tabcaption{\def\@captype{table}\caption}
\makeatother

\newtheorem{acknowledgement}[theorem]{Acknowledgement}

\newtheorem{axiom}[theorem]{Axiom}
\newtheorem{case}[theorem]{Case}
\newtheorem{conclusion}[theorem]{Conclusion}
\newtheorem{condition}[theorem]{Condition}
\newtheorem{conjecture}[theorem]{Conjecture}
\newtheorem{criterion}[theorem]{Criterion}
\newtheorem{example}[theorem]{Example}
\newtheorem{exercise}[theorem]{Exercise}
\newtheorem{notation}{Notation}
\newtheorem{solution}[theorem]{Solution}
\newtheorem{summary}[theorem]{Summary}

\newenvironment{proof}{\noindent {\bf
Proof.}}{\rule{3mm}{3mm}\par\medskip}
\newcommand{\remark}{\medskip\par\noindent {\bf Remark.~~}}
\newcommand{\pp}{{\it p.}}
\newcommand{\de}{\em}
\newcommand{\mad}{\rm mad}
\newcommand{\qf}{Q({\cal F},s)}
\newcommand{\qff}{Q({\cal F}',s)}
\newcommand{\qfff}{Q({\cal F}'',s)}
\newcommand{\f}{{\cal F}}
\newcommand{\ff}{{\cal F}'}
\newcommand{\fff}{{\cal F}''}
\newcommand{\fs}{{\cal F},s}
\newcommand{\s}{\mathcal{S}}
\newcommand{\G}{\Gamma}
\newcommand{\g}{(G_3, L_{f_3})}
\newcommand{\wrt}{with respect to }
\newcommand {\nk}{ Nim$_{\rm{k}} $  }

\newcommand{\gs}{(G, \sigma)}
\newcommand{\xeta}{{\mathtt{x}}^{\eta}}
\newcommand{\nm}{non-vanishing monomial }
\newcommand{\nicem}{nice monomial }

\newcommand{\brac}[1]{{\left(#1\right)}}
\newcommand{\sbrac}[1]{{\left[#1\right]}}
\newcommand{\set}[1]{\left\{#1\right\}}
\newcommand{\norm}[1]{{\left|#1\right|}}
\newcommand{\floor}[1]{{\left\lfloor #1 \right\rfloor}}
\newcommand{\ceil}[1]{{\left\lceil #1 \right\rceil}}

\newcommand{\q}{\uppercase\expandafter{\romannumeral1}}
\newcommand{\qq}{\uppercase\expandafter{\romannumeral2}}
\newcommand{\qqq}{\uppercase\expandafter{\romannumeral3}}
\newcommand{\qqqq}{\uppercase\expandafter{\romannumeral4}}
\newcommand{\qqqqq}{\uppercase\expandafter{\romannumeral5}}
\newcommand{\qqqqqq}{\uppercase\expandafter{\romannumeral6}}

\newcommand{\qed}{\hfill\rule{0.5em}{0.809em}}

\newcommand{\var}{\vartriangle}

\title{{\large \bf The Alon-Tarsi number of a  planar graph minus a matching  }}

\author{Jaros\l aw Grytczuk \thanks{Faculty of Mathematics and Information Science, Warsaw University of Technology, 00-662 Warsaw, Poland. E-mail: j.grytczuk@mini.pw.edu.pl. Supported by Polish National Science Center, Grant Number: NCN 2015/17/B/ST1/02660.} \and Xuding Zhu\thanks{Department of Mathematics, Zhejiang Normal University,  China.  E-mail: xudingzhu@gmail.com. Grant Number: NSFC 11571319.}}

\maketitle

\begin{abstract}
	
	This paper proves that every planar graph $G$ contains a matching $M$ such that the Alon-Tarsi number of $G-M$ is at most $4$. As a consequence, $G-M$ is $4$-paintable, and hence $G$ itself is $1$-defective $4$-paintable. This improves a result of Cushing and Kierstead [Planar Graphs are 1-relaxed, 4-choosable, 
		{\em European Journal of Combinatorics}  31(2010),1385-1397], who proved that every planar graph is $1$-defective $4$-choosable.


\end{abstract}

\section{Introduction}

In this paper we study the Alon-Tarsi number of special subgraphs of planar graphs, which is motivated by list coloring problems for planar graphs. A \emph{$k$-list assignment} of a graph $G$ is a mapping $L$ which assigns to each vertex $v$ of $G$ a set $L(v)$ of $k$ permissible colors. Given a $k$-list assignment $L$ of $G$, an \emph{$L$-colouring} of $G$ is a mapping $\phi$ which assigns to each vertex $v$ a colour $\phi(v) \in L(v)$ such that $\phi(u) \neq \phi(v)$ for every edge $e=uv$ of $G$. A graph $G$ is {\em $k$-choosable} if $G$ has an $L$-colouring for every $k$-list assignment $L$. The {\em choice number} of a graph $G$ is defined as
$$\ch(G) = \min\{k: \text{$G$ is $k$-choosable}\}.$$

Thomassen \cite{Thomassen94} proved that every planar graph $G$ is $5$-choosable. This bound is best possible, as proved by Voigt \cite{Voigt93}, who constructed the first non-$4$-choosable planar graph. Other examples of planar graphs with $\ch(G)=5$ can be found in \cite{ChoiKwon, Mirzakhani,Zhu2018}.

A natural question is ``how far" a planar graph can be from being $4$-choosable? One way to measure such distance is to consider defective list colouring, defined as follows.
A \emph{$d$-defective colouring} of a graph $G$ is a colouring of the vertices of $G$ such that each colour class induces a subgraph of maximum degree at most $d$. Thus, a $0$-defective colouring of $G$ is simply a proper colouring of $G$, while in a $1$-defective colouring a matching is allowed as a set of non-properly coloured edges.

Defective colouring of graphs was first studied by Cowen, Cowen and Woodall in \cite{CowenCW86}. They proved that every outerplanar graph is $2$-defective $2$-colourable, and that every planar graph is $2$-defective $3$-colourable. They also found examples of an outerplanar graph that is not $1$-defective $2$-colourable, a planar graph that is not $1$-defective $3$-colourable, and for every $d\geq 2$, a planar graph that is not $d$-defective $2$-colourable.

In a natural analogy to list colouring we may define defective list colouring of graphs. Given a $k$-list assignment $L$ of $G$, a \emph{$d$-defective $L$-colouring} of $G$ is a $d$-defective colouring $c$ of $G$ with $c(v) \in L(v)$ for every vertex $v$ of $G$.
A graph $G$ is \emph{$d$-defective $k$-choosable} if for any $k$-list assignment $L$ of $G$, there exists a $d$-defective $L$-colouring of $G$. Clearly, every $d$-defective $k$-choosable graph is $d$-defective $k$-colourable, however, the converse is not true. Nevertheless, the above mentioned results on defective colouring of planar graphs can be extended to defective list colouring. Eaton and Hull~\cite{EatonH99}, and \v{S}krekovski~\cite{Skrekovski99} proved independently that every planar graph is $2$-defective $3$-choosable, and every outerplanar graph is $2$-defective $2$-choosable. They both asked the natural question - whether every planar graph is $1$-defective $4$-choosable.
This problem is much more difficult. Only one decade later Cushing and Kierstead~\cite{CushingK10} answered this question in the affirmative, and their proof is rather complicated.

Another way to measure distance of a graph $G$ from being $k$-choosable is to consider the maximum degree of a subgraph that must be removed from $G$ in order to get a $k$-choosable graph. For example, we can ask the following question.
\begin{question} Is it true that every planar graph $G$ has a matching $M$ such that $G-M$ is $4$-choosable?
\end{question}

A positive answer to this question implies a stronger property than $1$-defective $4$-choosability of planar graphs. Indeed, we can specify a matching of non-properly coloured edges in advance, independently of the list assignment. Cushing and Kierstead \cite{CushingK10} proved that for any $4$-list assignment $L$ of $G$, there is a matching $M$ of $G$ such that $G-M$ is $L$-colourable. The matching $M$ constructed in the proof depends on the list assignment $L$.

Another way to extend the result of Cushing and Kierstead is to consider \emph{on-line} version of list colouring of graphs, defined through the following two person game \cite{Schauz09,Zhu09}. A \emph{$d$-defective $k$-painting game} on a graph $G$ is played by two players: Lister and Painter. Initially, each vertex is uncoloured and has $k$ tokens. In each round, Lister marks a chosen set $A$ of uncoloured vertices and removes one token from each marked vertex.
In response, Painter colours vertices in a subset $X$ of $A$ which induce a subgraph $G\sbrac{X}$ of maximum degree at most $d$. Lister wins if at the end of some round there is an uncoloured vertex with no more tokens left. Otherwise, after some round, all vertices are coloured, and then Painter wins the game. We say that $G$ is \emph{$d$-defective $k$-paintable} if Painter has a winning strategy in this game. A $0$-defective $k$-painting game is simply called a {\em $k$-painting game}, and a $0$-defective $k$-paintable graph is simply called {\em $k$-paintable}. The {\em paint number} of $G$ is defined as

$$\chi_{\p}(G) = \min\{k: \text{$G$ is $k$-paintable}\}.$$

It is not hard to see that every $d$-defective $k$-paintable graph is $d$-defective $k$-choosable. Indeed, let $L$ be some $k$-list assignment of $G$ with colours in the set $\sbrac{n}$. Suppose that Lister is playing in the following way. In the $i$-th round, for $i \in \sbrac{n}$, Lister marks the set $A_i=\set{v: i \in L(v), v \notin X_1\cup \ldots \cup X_{i-1}}$, where $X_j$ is the set of vertices coloured by Painter in the $j$-th round. But Painter has a winning strategy, so, he will eventually obtain a $d$-defective $L$-colouring of $G$. The converse is however not true. Though all planar graphs are $2$-defective $3$-choosable, an example of non-$2$-defective $3$-paintable planar graph has been constructed in \cite{GHKZ}. 

On the other hand, it is known that every planar graph is $5$-paintable \cite{Schauz09}, $3$-defective $3$-paintable \cite{GHKZ}, and $2$-defective $4$-paintable \cite{HanZ16}. For defective paintability of the family of planar graphs, Question \ref{q2} below is the only question remained open.

\begin{question}
	\label{q2}
Is it true that every planar graph is $1$-defective $4$-paintable?
\end{question}

More ambitiously, in analogy to Question 1, we may ask a similar question for paintability of planar graphs.

\begin{question}
Is it true that every planar graph $G$ has a matching $M$ such that $G-M$ is $4$-paintable?
\end{question}

  Our main result in this paper implies a positive answer to all of the stated questions. To formulate the main result, we need the following definitions. We associate to each vertex $v$ of $G$ a variable $x_v$. The \emph{graph polynomial} $P_G(\mathtt x)$ of $G$ is defined as
$$P_G(\mathtt{x}) = \prod_{uv\in E(G), u < v}(x_v-x_u),$$
where $\mathtt{x} = \{ x_v: v \in V(G)\}$ denotes the sequence of variables ordered accordingly to some fixed linear ordering $``<"$ of the vertices of $G$.
It is easy to see that a mapping $\phi: V \to \Bbb R$ is a proper colouring of $G$ if and only if
$P_G(\phi) \ne 0$, where $P_G(\phi)$ means to evaluate the polynomial at $x_v=\phi(v )$ for $v \in V(G)$.
Thus to find a proper colouring of $G$ is equivalent to find an assignment of $\mathtt{x}$ so that the polynomial evaluated at this assignment is non-zero.

Assume now that $P(\mathtt{x})$ is any real polynomial with variable set $X$. Let $\eta$ be a mapping which assigns to each variable $x$ a non-negative integer $\eta(x)$. We denote by ${\mathtt{x}}^{\eta}$ the monomial $\prod_{x \in X}  x^{\eta(x )}$ determined by mapping $\eta$, which we call then the \emph{exponent} of that monomial. Let $c_{P, \eta}$ denote the coefficient of ${\mathtt{x}}^{\eta}$ in the expansion of $P(\mathtt{x})$ into the sum of monomials. The celebrated Combinatorial Nullstellensatz of Alon \cite{AlonCN99} asserts that if $\sum_{x \in X}\eta(x)=\deg P(\mathtt{x})$ and $c_{P, \eta} \ne 0$, then for arbitrary sets $A_x$ assigned to variables $x \in X$, each consisting of $\eta(x)+1$ real numbers, there exists a mapping $\phi: X \to \Bbb R$ such that $\phi(x) \in A_x$ for each $x \in X$ and $P(\phi) \ne 0$.

Notice that a graph polynomial $P_G$ is \emph{homogenous}, which means that the exponents of each non-vanishing monomial sum up to the same value, which is equal the number of edges of $G$. Hence, condition $\sum_{x \in X}\eta(x)=\deg P_G(\mathtt{x})$ is satisfied by every non-vanishing monomial in $P_G$. In particular, Combinatorial Nullstellensatz implies that if $c_{P_G, \eta} \ne 0$ and $\eta(x_v) < k$ for all $v \in V$, then $G$ is $k$-choosable. This result was proved earlier by Alon and Tarsi \cite{AlonTarsi92}, who applied it, for instance, to demonstrate that planar bipartite graphs are $3$-choosable. It was then strengthened by Schauz \cite{Schauz09}, who showed that under the same assumptions, a graph $G$ is also $k$-paintable. 
Motivated by the above relations between list colourings and graph polynomials, Jensen and Toft \cite{JT1995} defined the {\em Alon-Tarsi number} $\AT(G)$ of a graph $G$ as $$\AT(G) = \min \{k: \text{$c_{P_G, \eta} \ne 0$ for some exponent $\eta$  with $\eta(x_v) < k$ for all $v \in V(G)$}\}.$$
As observed in \cite{Hefetz}, $\AT(G)$ has some distinct features, and it is of interest to study $\AT(G)$ as a separate graph invariant. Summarizing, for any graph $G$, we have $$\ch(G) \le \chi_{\p}(G) \leq \AT(G).$$The  gaps between these three parameters can be arbitrarily large. However, upper bounds for the choice number of many natural classes of graphs are also upper bounds for their Alon-Tarsi number. For example, as a strengthening of the result of Thomassen and the result of Schauz, it was shown in \cite{Zhu18} that every planar graph $G$ satisfies $\AT(G)\leq 5$.

In this paper we prove that every planar graph $G$ contains a matching $M$ such that $\AT(G-M)\leq 4$, which implies a positive answer to all three questions we formulated above.

\section{The main result}

The main result of this paper is the following theorem.

\begin{theorem}
	\label{thm-main0} Every planar graph $G$ contains a matching $M$ such that $\AT(G-M) \le 4$.
\end{theorem}
Before proceeding to the proof we need to fix some notation and terminology. For simplicity,   we write $c_{G,\eta}$ for $c_{P_G, \eta}$. We will also say that $\xeta$ is a monomial of a graph $G$ while formally it is a monomial in the graph polynomial $P_G$, and  the exponent $\eta$ assigns to each vertex $v$   a non-negative integer while formally the integer is assigned to $x_v$.  A monomial $\xeta$ of $G$ is {\em non-vanishing} if $c_{G,\eta} \ne 0$. By $d_H(v)$ we denote the degree of a vertex $v$ in a graph $H$. We will need the following definitions.

\begin{definition}
	\label{def-nice} Assume $G$ is a plane graph, $e=v_1v_2$ is a boundary edge of $G$, and $M$ is a matching in $G$. A monomial $\xeta$ of $G-e-M$ is {\em nice} for $(G,e,M)$ if the following conditions hold:
	\begin{enumerate}
		\item $\xeta$ is non-vanishing.
		\item $\eta(v_1)=\eta(v_2) =0$.
		\item    $\eta(v) \le 2 - d_M(v)$ for every other boundary vertex $v$.
		\item  $\eta(v) \le 3$ for each interior vertex $v$.
		\end{enumerate}	
\end{definition}
Notice that $d_M(v)=1$ if $v$ is covered by $M$, and $d_M(v)=0$ otherwise.
\begin{definition}
	Assume $G$ is a plane graph and $e=v_1v_2$ is a boundary edge of $G$. A matching $M$ of $G$ is \emph{valid} for $(G,e)$ if none of $v_1$ or $v_2$ is covered by $M$.
\end{definition}
If $\xeta$ is a \nicem for $(G,e, M)$, then let $\eta'(x)=\eta(x)$ except that $\eta'(x_{v_1})=1$.
As $P_G(\mathtt{x}) = (x_{v_1}-x_{v_2})P_{G-e}(\mathtt{x})$ and $\eta'(v_2)=0$, we know that $c_{G,\eta'} =c_{G-e, \eta} \ne 0$. Note that $\eta'(x_{v}) \leq 3$ for each vertex $v$. Thus Theorem \ref{thm-main0} follows from Theorem \ref{thm-2} below.
\begin{theorem}
	\label{thm-2}
	Assume $G$ is a plane graph and $e=v_1v_2$ is a boundary edge of $G$. Then $(G,e)$ has a valid matching $M$ such that there exists a \nicem $\xeta$ for $(G,e, M)$.
\end{theorem}

A variable $x$ is {\em dummy} in $P({\mathtt{x}})$ if  $x$ does not really occur in $P({\mathtt{x}})$, or equivalently, $\eta(x)=0$ for each non-vanishing monomial ${\mathtt{x}}^{\eta}$ in the expansion of $P$.    We shall frequently need to consider the summation and the product of polynomials.
By introducing dummy variables, we assume   the involved polynomials in the sum or the product have the same set of variables. For example,
we may view $x_2^2$ be the same as $x_1^0x_2^2x_3^0\ldots x_n^0$, that is,  $x_2^2={\mathtt{x}}^{\eta}$, where the variable set is $X=\{x_1, x_2, \ldots, x_n\}$, $\eta(x_2)=2$, and $\eta(x_i)=0$ for $i \ne 2$. We denote by $X$ the set of variables for   polynomials in concern. For two mappings $\eta_1, \eta_2$,
we write $\eta_1 \le \eta_2$ if $\eta_1(x) \le \eta_2(x)$ for all $x \in X$, and $\eta= \eta_2-\eta_1$ (respectively, $\eta= \eta_2+ \eta_1$) means  that $\eta(x) = \eta_2(x)-\eta_1(x)$ (respectively, $\eta(x) = \eta_2(x)+\eta_1(x)$) for all $x \in X$.

The following lemma just collects some simple observations needed for future  reference. Its proof amounts to routine checking and is therefore omitted.
\begin{lemma} The following properties hold for real polynomials and their monomial coefficients.
	\label{obs}
	\mbox{}	
	\begin{enumerate}
		\item If $P(\mathtt{x})=\alpha P_1(\mathtt{x})+ \beta P_2(\mathtt{x})$, then  $c_{P, \eta}=\alpha c_{P_1, \eta} + \beta c_{P_2, \eta}.$
		\item If $P({\mathtt{x}}) =  {\mathtt{x}}^{\eta'}P_1({\mathtt{x}})$, then
		$c_{P, \eta}=  c_{P_1, \eta-\eta'}.$
		\item If $P({\mathtt{x}}) =   {\mathtt{x}}^{\eta'}P_1({\mathtt{x}})$ and $\eta' \not\le \eta$, then $c_{P, \eta}=0$.
		\item  If $P({\mathtt{x}}) =  P_1({\mathtt{x}})P_2({\mathtt{x}})$ and for any $\eta'$ with $c_{P_2,\eta'}\ne 0$,  there is a dummy variable $x$ of $P_1({\mathtt{x}})$ such that $\eta'(x) \ne \eta(x)$,  then $c_{P, \eta}=0$.
		\item If $G$ is a graph and $c_{ G, \eta} \ne 0$, then $\sum_{x \in X} \eta(x) = |E(G)|$.
	\end{enumerate}
\end{lemma}

\noindent
{\bf Proof of Theorem \ref{thm-2}}
Assume the theorem is not true and $G$ is a minimum counterexample. It is not difficult to check that $G$ has at least $4$ vertices and is connected.
It is easy to see that if $(G,e,M)$ has a nice monomial and $G'$ is obtained from $G$ by deleting an edge, then $(G',e,M)$ also has a nice monomial. Thus we may assume that the boundary of $G$ is simple cycle.  

\bigskip
\noindent
{\bf Case 1 ($G$ has a chord.)} Let $f=xy$ be a chord in $G$. Let $G_1, G_2$ be the two $f$-components, that is, $G_1, G_2$ are induced subgraphs of $G$, where $V(G) = V(G_1) \cup V(G_2)$ and $V(G_1) \cap V(G_2) = \{x,y\}$), with $e \in G_1$. By the minimality of $G$,   $(G_1, e)$ has a valid   matching $M_1$ such that $(G_1,e, M_1)$ has a  \nicem $\mathtt{x}^{\eta_1}$, and $(G_2, f)$ has a valid matching $M_2$ such that $(G_2,f, M_2)$ has a \nicem $\mathtt{x}^{\eta_2}$.
Let $M=M_1 \cup M_2$. It is easy to see that $M$ is a valid matching in $(G, e)$.
Let $\eta=\eta_1+\eta_2$. We shall show that $\xeta$ is a \nicem for $(G,e,M)$. It is obvious that conditions (2)-(4) of a Definition \ref{def-nice} are satisfied by $\eta$. It remains to show that $c_{G-e-M, \eta} \ne 0$. It suffices to show that $c_{G-e-M, \eta} = c_{G_1-e-M_1, \eta_1} c_{G_2-f-M_2, \eta_2}$.

Assume $\mathtt{x}^{\eta'_1}$ is a non-vanishing monomial of $G_1-e-M_1$ and   $\mathtt{x}^{\eta'_2}$ is a non-vanishing monomial of $G_1-f-M_2$. Let $\eta'=\eta'_1+\eta'_2$.
	We shall show that $\eta' =\eta$ only if $\eta'_1=\eta_1$ and $ \eta'_2=\eta_2$. Assume $\eta' =\eta$. Since $\eta'_1(v) = 0$ for each vertex $v \in V(G_2)-\{x,y\}$, we have $\eta'_2(v) = \eta'(v) = \eta(v) =\eta_2(v)$ for each $v \in V(G_2)-\{x, y\}$. 
	As $\sum_{v \in V(G_2)} \eta'_2(v) = \sum_{v \in V(G_2)} \eta_2(v)$,
	we conclude that $\eta'_2(x) = \eta'_2(y)=0$. Hence $\eta'_2 = \eta_2$. This implies that $\eta'_1 = \eta- \eta_2 = \eta_1$.  This completes the proof of Case 1.

\bigskip
\noindent
{\bf Case 2 ($G$ has no chord.)} Assume $G$ has no chord and let $B(G)=(v_1, v_2, \ldots, v_n)$ be the vertices of the boundary of $G$ in a cyclic order. Let $G'=G-v_n$. By the minimality of $G$, $G'$ has a    valid matching $M$ such that $(G',e,M)$ has a \nicem $\mathtt{x}^{\eta'}$.

Let $v_1, u_1, u_2, \ldots, u_k, v_{n-1}$ be the neighbours of $v_n$. Since $G$ has no chord, each $u_i$ is an interior vertex of $G$. Let $$S({\mathtt{x}}) = (x_{v_n}-x_{v_1})(x_{v_{n-1}}-x_{v_n })( x_{u_1}-x_{v_n }) \ldots ( x_{u_k}-x_{v_n }).$$ Then
$P_{G-e-M}(\mathtt{x})=S({\mathtt{x}})P_{G'-e-M}(\mathtt{x})$.

Assume first that $n=3$. Let then
$\eta(x_v)=\eta'(x_v)$ for $v \notin \{u_1,u_2,\ldots, u_k\}$, $\eta(x_v) = \eta'(x_v)+1$ for
$v \in \{u_1,u_2,\ldots, u_k\}$, and $\eta(x_{v_n})=2$.
Let $\eta''(x_{v_n})=2$, $\eta''(x_{u_i})=1$ for $i=1,2,\ldots,k$, and $\eta''(x)=0$ for other $x$.
Then we may write $$S({\mathtt{x}}) = - {\mathtt{x}}^{\eta''} +x_{v_1}A(\mathtt{x}) + x_{v_2}B(\mathtt{x})+x_{v_n}^3C({\mathtt{x}})$$
for some polynimals $A(\mathtt{x})$,  $B(\mathtt{x})$ and $C({\mathtt{x}})$.  Let \begin{eqnarray*}
 P_1({\mathtt{x}}) &=& x_{v_1}A(\mathtt{x})P_{G'-e-M}({\mathtt{x}}),\\
 P_2({\mathtt{x}}) &=& x_{v_2}B(\mathtt{x})P_{G'-e-M}({\mathtt{x}}),\\
 P_3({\mathtt{x}}) &=& x^3_{v_n}C(\mathtt{x})P_{G'-e-M}({\mathtt{x}}).
\end{eqnarray*}
As $\eta(x_{v_1})=\eta(x_{v_2})=0$ and $\eta(x_{v_n})=2$, it follows from (3) of Lemma \ref{obs} that $c_{P_1,\eta}=c_{P_2,\eta}=c_{P_3,\eta}=0$. By (1) and (2) of Lemma \ref{obs},
$c_{G-e, \eta} = -c_{G'-e, \eta'} \ne 0$. Hence $\eta$ is nice for $(G,e. M)$.

Assume now that $n \ge 4$. We say a monomial $\mathtt{x}^{\tau}$ for $G'-e-M$  is {\em special} if
\begin{itemize}
	\item  $\tau(v_1) = \tau(v_2)=0$.
	\item $\tau(v_{n-1}) \le 1-d_M(v_{n-1})$.
	\item $\tau(v) \le 2-d_M(v)$ for every other boundary vertex $v$, except that there may be at most one index $i \in \{1,2,\ldots,k\}$ such that
	$\tau(u_i) = 3-d_{M}(u_i)$.
	\item $\tau(v) \le 3$ for each interior vertex  $v$.
\end{itemize}

\bigskip
\noindent
{\bf Subcase 2(i) (There is a non-vanishing special monomial in $G'-e-M$.)} Assume that $c_{ {G'-e-M}, \tau}  \ne 0$ for some special monomial $\mathtt{x}^{\tau}$ in $G'-e-M$. For $i \in \{1,2,\ldots, k\}$, we say $u_i$ is {\em saturated} if $\tau(u_i) = 3$. By the definition of special monomial, we know that there is at most one  $u_i  $ that is saturated. Moreover, if $u_i$ is saturated, then $d_{M}(u_i) =0$. Let
\[
M'= \begin{cases} M, & \text{ if no $u_i$ is saturated}, \cr
M \cup \{u_iv_n\}, & \text{ if $u_i$ is saturated}. \cr
\end{cases}
\]
It follows from the definition that $M'$ is a valid matching in $(G, e)$. Let
\[
\eta(v) = \begin{cases} \tau(v), &\text{if $v \notin \{u_1, u_2, \ldots, u_k, v_{n-1}\}$ or $v=u_i$   is   satuarated } \cr
\tau(v)+1, &\text{ if $v \in \{u_1, u_2, \ldots, u_k, v_{n-1}\}$  is not saturated}, \cr
1, &\text{if $v = v_n$.}
\end{cases}
\]
Let $\tau'$ be a mapping defined as $\tau'( v_n ) = \tau'( v_{n-1} )=1$, $\tau'(u_j)=1$ if $u_j$ is not saturated, and $\tau'(v)=0$ for other vertices $v$. Let $$\tilde{S}({\mathtt{x}}) = (x_{v_n}-x_{v_1})(x_{v_{n-1}}-x_{v_n }) \prod_{ \text{  $u_i$ is not   saturated}} ( x_{u_i}-x_{v_n }).$$
Then $$P_{G-e-M', \eta} = \tilde{S}({\mathtt{x}}) P_{G'-e-M, \tau}$$ and
$$\tilde{S}({\mathtt{x}}) =  {\mathtt{x}}^{\tau'} +x_{v_1}A(\mathtt{x}) + x_{v_n}^2 B(\mathtt{x})$$
for some polynomials $A(\mathtt{x})$ and $B(\mathtt{x})$.

Let $P({\mathtt{x}}) = {\mathtt{x}}^{\tau'}P_{G'-e-M}({\mathtt{x}})$, $P_1(\mathtt{x})=  x_{v_1}A(\mathtt{x})P_{G'-e-M}(\mathtt{x})$, and $P_2(\mathtt{x})= x^2_{v_n}B(\mathtt{x})P_{G'-e-M}(\mathtt{x}) $.
Then
$$P_{G-e-M'}({\mathtt{x}}) = P({\mathtt{x}})+P_1({\mathtt{x}})+P_2({\mathtt{x}}).$$ As $\eta(v_1)=0$ and $\eta(v_n)=1$, it follows from (3) of Lemma \ref{obs} that
$c_{P_1, \eta  } =  c_{P_2, \eta  }=0.$
By (1) and  (2) of Observation \ref{obs}, we have
$c_{  {G-e-M'}, \eta}=c_{ P, \eta} =   c_{  {G'-e-M}, \tau} \ne 0.$
Hence $\xeta$ is a \nicem for $(G,e,M')$.

\bigskip
\noindent
{\bf Subcase 2(ii) (There is no non-vanishing special monomial in $G'-e-M$.)} We assume now that $c_{ {G'-e-M}, \tau} = 0$ for every special monomial $\mathtt{x}^{\tau}$ of $G'-e-M$. Recall that $\mathtt{x}^{\eta'}$ is a  \nicem for $(G',e, M)$. Let  $\eta(x_v)=\eta'(x_v)$ for $v \notin \{u_1,u_2,\ldots, u_k\}$, $\eta(x_v) = \eta'(x_v)+1$ for
$v \in \{u_1,u_2,\ldots, u_k\}$, and $\eta(x_{v_n})=2$. We shall show that $\xeta$ is a \nicem for $G-e-M$.

It is obvious that Conditions (2)-(4) of Definition \ref{def-nice} are satisfied by $\xeta$. It remains to show that $\xeta$ is non-vanishing. Note that
$P_{G-e-M}(\mathtt{x})=S({\mathtt{x}})P_{G'-e-M}(\mathtt{x}).$ Let $\tau'(x_{v_n})=2$, $\tau'(x_{u_i})=1$ for $i=1,2,\ldots,k$, and $\tau'(x) = 0$ for other $x$. Then
$\eta = \tau'+\eta'$ and
$$S({\mathtt{x}})=  {\mathtt{x}}^{\tau'}  - \sum_{i=1}^k x_{v_n}^2x_{v_{n-1}}  \prod_{j\ne i}x_{u_j} +x_{v_1}A(\mathtt{x})+x_{v_n}^3B({\mathtt{x}})$$
for some polynimals $A(\mathtt{x})$ and $B(\mathtt{x})$. Let
\begin{eqnarray*}
 P(\mathtt{x}) &=& {\mathtt{x}}^{\tau'} P_{G'-e-M}  (\mathtt{x}), \\
 P_1(\mathtt{x}) &=&  \sum_{i=1}^k x_{v_n}^2x_{v_{n-1}}  \prod_{j\ne i}x_{u_j}  P_{G'-e-M}(\mathtt{x}),\\
 P_2(\mathtt{x}) &=& x_{v_1}A(\mathtt{x}) P_{G'-e-M}(\mathtt{x}),\\   P_3(\mathtt{x}) &=& x_{v_n}^3B(\mathtt{x}) P_{G'-e-M}(\mathtt{x}).
\end{eqnarray*}
As $\eta(x_{v_1})=0$ and $\eta(x_{v_n})=2$, it follows from (3) of Lemma \ref{obs} that $c_{ P_2, \eta  }  =c_{P_3, \eta}=0.$

For $i=1,2,\ldots,k$, let
$$P_{1,i} =x_{v_n}^2x_{v_{n-1}}  \prod_{j\ne i}x_{u_j}  P_{G'-e-M}(\mathtt{x})$$
and let
\[
\tau_i(v) = \begin{cases} \eta'(v)-1, &\text{if $v=v_{n-1}$}, \cr
\eta'(v)+1, &\text{if $v=u_i$}, \cr
\eta'(v), &\text{ otherwise}.
\end{cases}
\]
Then
$$c_{P_{1,i}, \eta} =  -    c_{  {G'-e-M}, \tau_i} $$ and hence
$$c_{P_1, \eta} =  - \sum_{i=1}^k  c_{  {G'-e-M}, \tau_i}.$$
For each $i \in \{1,2,\ldots,k\}$,  $\mathtt{x}^{\tau_i}$ is a special monomial for $G'-e-M$. Hence $ c_{  {G'-e-M}, \tau_i} = 0$, and consequently,
$c_{P_1, \eta } = 0$. As
$$P_{G-e-M}({\mathtt{x}}) = P({\mathtt{x}}) + P_1({\mathtt{x}})+P_2({\mathtt{x}})+P_3({\mathtt{x}}),$$ by (1) and (2) of Lemma \ref{obs} we get $c_{ {G-e-M}, \eta}=c_{P, \eta} =   c_{ {G'-e-M}, \eta'} \ne 0$. This finishes the proof of Case 2, and completes the proof of the theorem.
\qed

 \section{Some remarks}

 Assume $f: V(G) \to \{1,2,\ldots, \}$ is a function which assigns to each vertex $v$ of $G$ a positive integer. We say $G$ is {\em $f$-choosable} if for any list assignment $L$ with $|L(v)| = f(v)$ for every vertex $v$, $G$ is $L$-colourable. The {\em $f$-painting game } is defined in the same way as the $k$-painting game, except that initially,  instead of $k$ tokens, each vertex $v$ has $f(v)$ tokens. We say $G$ is {\em $f$-paintable} if Painter has a winning strategy in the $f$-painting game on $G$. We say $G$
  {\em $f$-Alon-Tarsi}, or {\em $f$-$\AT$} for short, if $P_G$ has a non-vanishing monomial $\xeta$ with $\eta(v) < f(v)$ for each vertex $v$.

In the proof of Theorem \ref{thm-2}, in Case 1, instead of adding the edge $u_iv_n$ to the matching $M'$, we may increase the power of $x_{u_i}$ by $1$. Then the resulting monomial  is non-vanishing in $P_{G-e}$. Thus a slight modification of the proof of Theorem \ref{thm-2} proves the following theorem, which is a strengthening of the result that every planar graph $G$ has $\AT(G) \le 5$.

\begin{theorem}
	\label{thm-main2}
	Assume $G$ is a planar graph. Then $G$ has a matching $M=\{(x_i, y_i): i=1,2,\ldots, p\}$ (note that edges in $M$ are oriented) such that $G$ is $f$-$\AT$, where   $f: V(G) \to \{4,5\}$ is defined as $f(x_i)=5$ and $f(v)=4$ for $v \in V-\{x_1,x_2,\ldots, x_p\}$. Consequently, $G$ is $f$-paintable, and hence $f$-choosable.
\end{theorem}
Note that $|V(G)| > 2p$. Thus we have the following corollary, which is the strengthening of the $5$-choosability of planar graphs.

\begin{corollary}
	\label{cor3}
	Every planar graph $G$ has a subset $X$ of vertices with $|X| < |V(G)|/2$ such that if $L$ is a list assignment which assigns to each vertex in $X$ five permissible colours and assigns to each other vertex  four permissible colours, then $G$ is $L$-colourable.
\end{corollary}

A {\em signed graph} is a pair $(G, \sigma)$, where $G$ is a graph and $\sigma$ is a {\em signature} of $G$  which assigns to each edge $e$ of $G$ a sign $\sigma_e \in \{1,-1\}$. A {\em proper colouring} of $(G, \sigma)$ is a mapping $f$ which assigns to each vertex $v$ an integer $f(v)$ so that for each edge $e=xy$, $f(x) \ne \sigma_e f(y)$. The chromatic number $ \chi \gs $ of $\gs$ is the minimum integer $k$ such that for any  set $S$ of $k$ integers, $\gs$ has a proper colouring using colours from $S$.
The choice number of $(G, \sigma)$ is the minimum integer $k$ such that for every   $k$-list assignment $L$ of $G$, there is a proper $L$-colouring of  $(G, \sigma)$. The polynomial associated to $(G, \sigma)$ is defined as $$P_{G, \sigma}(\mathtt{x}) = \prod_{ u \sim v, u < v}(x_v- \sigma(uv)x_u).$$
The Alon-Tarsi number of signed graphs is defined similarly. Then all the arguments in the proof of Theorem \ref{thm-2} works. Hence we have the following result.

\begin{theorem}
	\label{thm-main3} If $(G, \sigma)$ is a signed  planar graph, then $G$ has a matching $M$ such that  $AT(G-M, \sigma) \le 4$. Consequently, $(G-M, \sigma)$ is $4$-choosable, and $(G, \sigma)$ itself is $1$-defective $4$-choosable.
\end{theorem}
Corollary \ref{cor3} also works for signed planar graphs.
\begin{corollary}
		\label{cor4}
	Every signed planar graph $(G, \sigma)$ has a subset $X$ of vertices with $|X| < |V(G)|/2$ such that if $L$ is a list assignment which assigns to each vertex in $X$ five permissible colours and asigns to each other vertex  four permissible colours, then $(G, \sigma)$ is $L$-colourable.
\end{corollary}

\bigskip
\noindent
{\bf Acknowledgement}
  Jaros\l aw Grytczuk 
would like to thank Zhejiang Normal University for great hospitality, and to all people from the Center for Discrete Mathematics  for a wonderful atmosphere during his visit in Jinhua, where this research was carried out. His visit was supported by the 111 project of the Ministry of Education of China on ``Graphs and Network Optimization".

\bibliographystyle{plain}
\bibliography{paper1}

\begin{thebibliography}{10}

\bibitem{AlonCN99}
N.~Alon.
\newblock Combinatorial nullstellensatz.
\newblock {\em Combinatorics, Probability, and Computing}, 8:7--29, 1999.

\bibitem{AlonTarsi92}
N.~Alon and M.~Tarsi.
\newblock Colorings and orientations of graphs.
\newblock {\em Combinatorica}, 12:125--134, 1992.

\bibitem{ChoiKwon}
H.~Choi and Y.~S. Kwon.
\newblock On {$t$}-common list-colorings.
\newblock {\em Electron. J. Combin.}, 24(3):Paper 3.32, 10, 2017.

\bibitem{CowenCW86}
L.~J. Cowen, R.~H. Cowen, and D.~R. Woodall.
\newblock Defective colorings of graphs in surfaces: Partitions into subgraphs
  of bounded valency.
\newblock {\em Journal of Graph Theory}, 10(2):187--195, 1986.

\bibitem{CushingK10}
W.~Cushing and H.~A. Kierstead.
\newblock Planar graphs are 1-relaxed, 4-choosable.
\newblock {\em European Journal of Combinatorics}, 31(5):1385 -- 1397, 2010.

\bibitem{EatonH99}
N.~Eaton and T.~Hull.
\newblock Defective list colorings of planar graphs.
\newblock {\em Bulletin of the Institute of Combinatorics and its
  Applications}, 25:79--87, 1999.

\bibitem{GHKZ}
G.~Gutowski, M.~Han, T.~Krawczyk, and X.~Zhu.
\newblock Defective 3-paintability of planar graphs.
\newblock {\em Electron. J. Combin.}, 25(2):Paper 2.34, 20, 2018.

\bibitem{HanZ16}
M.~Han and X.~Zhu.
\newblock Locally planar graphs are 2-defective 4-paintable.
\newblock {\em European Journal of Combinatorics}, 54:35--50, 2016.

\bibitem{Hefetz}
D.~Hefetz.
\newblock On two generalizations of the alon-tarsi polynomial method.
\newblock {\em Journal of Combinatorial Theory Ser. B}, 101(2):403--414, 2011.

\bibitem{JT1995}
T.~Jensen and B.~Toft.
\newblock {\em Graph Coloring Problems}.
\newblock Wiley, New York, 1995.

\bibitem{Mirzakhani}
M.~Mirzakhani.
\newblock A small non-{$4$}-choosable planar graph.
\newblock {\em Bull. Inst. Combin. Appl.}, 17:15--18, 1996.

\bibitem{Schauz09}
U.~Schauz.
\newblock {M}r.\ {P}aint and {M}rs.\ {C}orrect.
\newblock {\em Electronic Journal of Combinatorics}, 16(1):R77:1--18, 2009.

\bibitem{Thomassen94}
C.~Thomassen.
\newblock Every planar graph is 5-choosable.
\newblock {\em Journal of Combinatorial Theory, Series B}, 62(1):180--181,
  1994.

\bibitem{Voigt93}
M.~Voigt.
\newblock List colourings of planar graphs.
\newblock {\em Discrete Mathematics}, 120(1):215--219, 1993.

\bibitem{Skrekovski99}
R.~\v{S}krekovski.
\newblock List improper colourings of planar graphs.
\newblock {\em Combinatorics, Probability and Computing}, 8(3):293--299, 1999.

\bibitem{Zhu18}
X.~Zhu.
\newblock Alon-{T}arsi number of planar graphs.
\newblock {\em Journal of Combinatorial Theory Ser. B},
  https://doi.org/10.1016/j.jctb.2018.06.00.

\bibitem{Zhu09}
X.~Zhu.
\newblock On-line list colouring of graphs.
\newblock {\em Electronic Journal of Combinatorics}, 16(1):R127:1--16, 2009.

\bibitem{Zhu2018}
X.~Zhu.
\newblock A refinement of choosability of graphs, 2018.

\end{thebibliography}

\end{document}